\documentclass[11pt]{article}
\usepackage[
english]{babel}
\usepackage{amsmath,amssymb,xcolor,theorem,bbm,times}

\usepackage{accents}
\usepackage{xcolor}
\definecolor{blue}{rgb}{0,0,1}
\usepackage[pagebackref,bookmarks]{hyperref} 
\hypersetup{pagebackref=true,bookmarks=true,
       hyperfigures=true,
	unicode=true,
	colorlinks=true,
	citecolor=blue,
	linkcolor=blue,
	anchorcolor=red 
}

\usepackage[textsize=footnotesize,color=blue!40, bordercolor=white]{todonotes}
\newcommand\todoi[1]{\todo[inline]{#1}}
\catcode`\<=\active \def<{
\fontencoding{T1}\selectfont\symbol{60}\fontencoding{\encodingdefault}}
\catcode`\>=\active \def>{
\fontencoding{T1}\selectfont\symbol{62}\fontencoding{\encodingdefault}}
\catcode`\|=\active \def|{
\fontencoding{T1}\selectfont\symbol{124}\fontencoding{\encodingdefault}}
\newcommand{\nobracket}{}
\newcommand{\nocomma}{}
\newcommand{\nosymbol}{}

\newcommand{\tmem}[1]{{\em #1\/}}
\newcommand{\tmop}[1]{\ensuremath{\operatorname{#1}}}

\newcommand{\tmtextit}[1]{{\itshape{#1}}}
\newtheorem{theorem}{Theorem}[section]

\newtheorem{corollary}[theorem]{Corollary}
\newtheorem{definition}[theorem]{Definition}

\newtheorem{lemma}[theorem]{Lemma}
\newtheorem{remark}[theorem]{Remark}

\usepackage{fancyhdr}

\newcommand{\cali}{\ensuremath{\mathcal{I}} }

\newcommand{\call}{ {\mathcal{L}}}
\newcommand{\calc}{C}

\newcommand{\nod}{{\noindent}}

\newcommand{\bu}{\ensuremath{\bullet } }
\newcommand{\nat}{\ensuremath{\mathbbm{N}}}
\newcommand{\re}{\ensuremath{\mathbbm{R}}}

\newcommand{\rest}{{\upharpoonright}}
\newcommand{\emp}{{\varnothing}}
\newcommand{\pa}[1]{\ensuremath{\langle #1 \rangle}}

\newcommand{\calt}{T}

\newcommand{\calu}{\ensuremath{\mathcal{U}}}

\newcommand{\vp}{\ensuremath{\varphi}}

\newcommand{\la}{{\langle}}
\newcommand{\ra}{{\rangle}}

\newcommand{\dom}{\text{dom}}

\newcommand{\power}{P}

\newcommand{\ie}{{\itshape{i.e.}}{\hspace{0.25em}}}

\newcommand{\etc}{{\itshape{etc.}}{\hspace{0.25em}}}
\newcommand{\via}{{\itshape{via }}{\hspace{0.25em}}}
\newcommand{\cf}{{\itshape{cf. }}{\hspace{0.25em}}}

\newcommand{\Equi}{{\,\Longleftrightarrow\,}}

\newcommand{\ex}{{\exists}}

\newcommand{\sset}{ {\, \subseteq  \,} }

\newcommand{\imp}{{\,\longrightarrow\,}}

\newcommand{\fin}[1]{\ensuremath{[  ]^{< \omega}}}

\newcommand{\pf}{{\noindent}{\textbf{Proof: }}}
\newcommand{\dfs}{\ensuremath{=_{\ensuremath{\operatorname{df}}}}}

\newcommand{\qed}{ \mbox{ } \hfill Q.E.D. }

\def\mid{\, | \, }

\def\m``{\mbox{``}}
\def\mq{\mbox{''}}

\def\le0{L[E_{0}]}

\def\ii{\iota}

\def\l{\lambda}

\def\colr{\color{red}}
\def\mb{\mbox{}}
\def\power{\mathcal{P}}
\def\re{\mathbbm{R}}
\def\calt{\cal{T}}
\def\pistol{\mathparagraph}

\def\WCL{\mbox{WCL}}
\def\ZFC{\mbox{ZFC}}

\begin{document}

\title{ { When cardinals determine the power set: inner models and H\"artig quantifier logic\thanks{The first author has received 
funding from the Academy of Finland (grant No 322795) and from 
the European Research Council (ERC) under the European Union’s Horizon 2020 
research and innovation programme (grant agreement No 101020762).}}}

\author{J. V\"a\"an\"anen\\
University of Helsinki\\
 and \\
University of Amsterdam \and P.D. Welch\\
University of Bristol}

\date{ \today }

\maketitle

\begin{abstract}
   
  We make use of some observations on the core model, for example assuming $V=L [ E ]$, that if there is no  inner model with a
  Woodin cardinal, and $M$  is an inner model with the same cardinals as
  $V$, then $V=M$. We conclude in this latter situation that ``$x=\mathcal{P}
  ( y )$'' is $\Sigma_{1} ( \tmop{Card} )$ where $\tmop{Card}$ is a predicate
  true of just the infinite cardinals. It is known that this implies the validities of second order logic are reducible to $V_I$ the set of validities of the H\"artig quantifier logic.

  We draw some further conclusions on the
  L\"owenheim number, $\ell_{I}$ of the latter  logic: that  if 
  no $L[E]$ model has a cardinal strong up to an $\aleph$-fixed point, and $\ell_{I}$
  is less than the least weakly inaccessible $\delta$, then (i) $\ell_I$ is a limit of measurable cardinals of $K$; (ii)
  the Weak Covering Lemma holds at $\delta$.
\end{abstract}

\section{Introduction}

The predicates ``$\alpha$ is a cardinal'', which we shorten $\tmop{Card}(\alpha)$, and ``$x$ is the power-set of $y$'', 
which we shorten $x=\mathcal{P}(y)$, are  paradigm examples of $\Pi_1$-predicates of set theory. In some models of set theory it may happen that the latter is $\Sigma_1$ in the former.
 In a sense, cardinal numbers then determine all power-sets, or in other words, the power-set operation is generalized recursive in the cardinal-predicate. We analyse this state of affairs and show that it is consistent with large cardinals below the first Woodin cardinal.

Let us use $\Sigma_1(\tmop{Card})$ to denote the class of $\Sigma_1$-formulas in the extended vocabulary $\{\in,\tmop{Card}\}$. Our concern in this paper is the proposition
\begin{equation}\label{star}
\mbox{$x=\mathcal{P}(y)$ is $\Sigma_1(\tmop{Card})$-definable}
\end{equation}
introduced in \cite{MR510714}. The meaning of (\ref{star}) is that for some $\Delta_0$-formula $\phi(z,x,y)$ of the vocabulary $\{\in,\tmop{Card}\}$ the equivalence
$$\forall x\forall y(x=\mathcal{P}(y)\leftrightarrow\exists z\phi(z,x,y))$$ holds.

Because $\mathcal{P}(\alpha)\cap L\subseteq L_{|\alpha|^+}$, it is clear that (\ref{star}) holds if $V=L$. It needs a little more to see that (\ref{star}) holds also if $V=L^\mu$. The purpose of this paper is to show that (\ref{star}) holds if $V=L[E]$ for a sequence $E$ of extenders and there is no inner model with a Woodin cardinal. This establishes the consistency of (\ref{star}) with large cardinals below the first Woodin cardinal, relative to the existence of such cardinals. It remains an open question how far this result can be extended. In particular, the question whether (\ref{star}) is consistent with a supercompact cardinal relative to the consistency of a supercompact cardinal, posed in \cite{MR694269}, remains open.

Note that the axiom (\ref{star}) by no means limits the size of the continuum. For example, (\ref{star}) is consistent with the negation of the Continuum Hypothesis, relative to the consistency of $ZF$ \cite{MR694269}. The axiom fails if we add a Cohen real \cite{MR510714}. For a stronger result, let $HC'$ be the set of sets of hereditary cardinality less than $2^\omega$. It is consistent relative to the consistency of $ZF$ that $$HC'\prec_{\Sigma_1(\tmop{Card})}V,$$ i.e. that $HC'$ reflects all $\Sigma_1(\tmop{Card})$-predicates, and then certainly  (\ref{star}) fails \cite{MR1928384}.


%
%


The origin of the hypothesis (\ref{star}) is in the model theory of generalized quantifiers. Let $I$
be the following generalized quantifier, known as the H\"artig-quantifier, or {\em equicardinality} quantifier:
$$Ixy\phi(x,\vec{a})\psi(y,\vec{b})\iff|\{x : \phi(x,\vec{a})\}|=\{y : \psi(y,\vec{b})\}|.$$
Let $L^I$ denote the extension of first order logic by the quantifier $I$. Early results on $L^I$ indicated that it is quite a strong logic and the question arose whether it is as strong as second order logic. This cannot be literally true for in a finite unary vocabulary the logic $L^I$ is decidable while second order logic is certainly not. However, there is a deeper sense in which the answer depends on (\ref{star}). To see what this means we need to introduce some notation.

Suppose $L^*$ is an abstract logic. We are mainly interested in the cases that $L^*$ is $L^I$ or  second order logic.
A class $K$ of models of a fixed vocabulary $\tau$ is \emph{$L^*$-definable} if there is an $L^*$-sentence $\phi$ such that $K$ is the class of  models of $\phi$. A class $K$ of models, again of a fixed vocabulary $\tau$, is $\Sigma(L^*)$  (-definable) if there is an $L^*$-sentence $\phi$, with a possibly larger vocabulary, such that $K$ is the class of relativized reducts of models of $\phi$. Finally, a class $K$ is said to be $\Delta(L^*)$ (-definable) if both $K$ and its complement in the class of all models of the vocabulary $\tau$ are $\Sigma(L^*)$. It was shown in \cite{MR457146} that we can regard $\Delta(L^*)$ as an abstract logic, as it is closed under finite unions and intersections as well as complementation, and (in a sense which is made precise in \cite{MR457146}) also under quantification. It is called the $\Delta$-{\em extension} of $L^*$. Intuitively, $\Delta(L^*)$ is the closure of $L^*$ under ``recursive'' operations. 

An important example, due to \cite{MR244012}, is the following: The class $K$ of models $(A,<)$ such that $<$ well-orders $A$, is $\Delta(L^I)$-definable. First of all the  complement of $K$ is clearly $\Sigma(L^I)$. To see that $K$ itself is $\Sigma(L^I)$ we consider the conjunction $\theta$ of the following $L^I$-sentences:
\begin{enumerate}
\item ``$<$ linearly orders the set $A$''
\item $\forall x\forall y\forall z(x<y\to (R(x,z)\to R(y,z)))$
\item $\forall x\forall y(x<y\to\neg IuvR(x,u)R(y,v))$.
\end{enumerate} It is easy to see that an ordered set $(A,<)$ is a well-order  if and only if it is a reduct of a model of $\theta$. This shows that $K$ is $\Delta(L^I)$.

Two logics are said to be {\em equivalent} if they have the same definable  model classes. For first order logic and $L_{\omega_1\omega}$, the logic and its $\Delta$-extension are equivalent---a consequence of the Craig Interpolation Theorem \cite{MR0406772}.  A logic is called  $\Delta$-closed if it is equivalent to its own $\Delta$-closure. Of course, the $\Delta$-extension of any logic is itself $\Delta$-closed. The logics are called $\Delta$-equivalent if their $\Delta$-extensions are equivalent. 

Historically the first example of $\Delta$-equivalence was the observation that the logic $L(Q_0)$, with the quantifier $$Q_0x\phi(x,\vec{s})\iff
|\{b : \phi(b,\vec{a})\}|\ge\aleph_0,$$ and weak second order logic $L^2_w$, with second order quantifiers over finite sets, relations and functions, are $\Delta$-equivalent. Moreover, the infinitary logic $L_{\mbox{\tiny HYP}}$ with conjunctions and disjunctions over recursive sets of formulas is $\Delta$-equivalent to $L(Q_0)$ and $L^2_w$.

The interest in the $\Delta$-operation stems from the fact that it preserves many model-theoretic properties. To see what this means, suppose $L^*$ and $L^{\mbox{\tiny +}}$ are $\Delta$-equivalent. Then 

\begin{enumerate}
\item $L^*$ satisfies the Compactness Theorem if and only if $L^{\mbox{\tiny +}}$ does. This generalizes to weaker form of compactness, such as $\kappa$-compactness.
\item  $L^*$ satisfies the Downward L\"owenheim-Skolem  Theorem\footnote{If a sentence has a model it has a countable model.} if and only if $L^{\mbox{\tiny +}}$ does. This generalizes to modifications  of the Downward L\"owenheim-Skolem  Theorem, such as the Downward L\"owenheim-Skolem  Theorem down\footnote{If a sentence has a model it has a model of cardinality $\leq\kappa$.} to $\kappa$.
\item $L^*$ is effectively axiomatizable\footnote{There is a G\"odel-numbering of all the formulas and the  set of G\"odel-numbers of valid formulas (called the {\em decision problem} of the logic) which is r.e.} if and only if $L^{\mbox{\tiny +}}$ is. More generally, the decision problems of $L^*$ and $L^{\mbox{\tiny +}}$ are recursively isomorphic.

\end{enumerate}

After this short introduction to the $\Delta$-operation we can state the connection between (\ref{star}) and the logic $L^I$: The following conditions are equivalent:
\begin{enumerate}
\item $x=\mathcal{P}(y)$ is $\Sigma_1(\tmop{Card})$-definable.
\item $L^I$ and second order logic are $\Delta$-equivalent.
\end{enumerate}




Our applications here concern semantic concepts related $L^I$ in particular to classifying the complexity of the validities set $V_{I}$ already mentioned and to the {\em L\"owenheim number} of $L^I$:
\begin{definition}
  The {\tmem{L\"owenheim number}} of a logic $\ell_{L^{\ast}}$ is the least
  cardinal $\kappa  $ so that for any $\varphi \in L^{\ast}$, if $\vp$ has a
  model then it has a model of cardinality $\leq   \kappa$.
\end{definition}

As intimated we apply this principally here to the logic $L^I$.  In this
case we write $\ell_I$ for this L\"owenheim number. A number of facts are easily established concerning the cardinality of $\ell_I$.

Note 1: $\ell_I$  is  of cofinality $\omega$, and is moreover a fixed point of the $\aleph$ function, thus $\ell_I =
\aleph_{\ell_I}$.

Note 2: This is not to be confused with the L\"owenheim-Skolem-Tarski number,
$\tmop{LST} ( L^{\ast} )$, which involves the notion of structures suitable
for a logic $L^{\ast}$ to have elementary substructures of size less than the
cardinal. This LST number may not exist (large cardinals are required for
example to show that it exists for second order logic $L^{2}$, or for $L^I$). We also have that $\tmop{LST}(L^I)$ is larger than
$\ell_I$, as can be easily seen by building a model $M$ from the countable number of witnesses of $\ell_I$ being what it is defined to be, and noting that $M$ cannot have an elementary submodel with respect to $LI$ of cardinality $<\ell_I$.   
Note 3: in \cite{MaVa11} at Theorem 21, it is shown that:\\

\nod{\bf Theorem}
  \label{T1.2}$\tmop{Con} ( \tmop{\ZFC}  + \exists \kappa (
  \kappa$ {\em  supercompact})) $\imp$
  
  $\tmop{Con} ( \tmop{\ZFC} + \tmop{LST} ( L ^I  ) =  \delta$, {\em {the}  
  {least}   {weakly}   {inaccessible}   {cardinal}} $)$.\\

\nod For the rest of this paper we let $\delta$ be the least weakly inaccessible
cardinal, if it exists.  If we write something such as ``$\theta < \delta$'' this is taken to assert that $\delta$ exists as well the truth of the inequality shown.  The following problem was open for several decades.\\

\nod {\bf Problem 2.5:} (\cite{MR654791}) {\tmem{If there are weakly inaccessible cardinals, can we have $\ell_I<\delta$?}}\\

As $\ell_I < \tmop{LST} ( L^I ) $ (see Note 2) we thus have by the last theorem the
consistency of $\ell_I < \delta$ relative to that of a supercompact.
That $\tmop{LST} ( {L^I} )$ exists, and is equal to $\delta$ is of
necessity a large cardinal notion: as is further shown in {\cite{MaVa11}} for
$\lambda \in \tmop{Card}$, $\lambda > \tmop{LST} ( {L^I} )$ we have the
failure of $\Box_{\lambda}$. The failure everywhere of this combinatorial
principle is known to imply large cardinals in inner models, and, as only an
example, projective determinacy.

This begs the question of a lower bound to the consistency strength of the simpler assertion $\ell_I<\delta$. Again in \cite{MR654791} it is shown there is no generic extension of $L$, or $L^{\mu}$ in
which $\ell_I < \delta$.
We note here that essentially what was shown there is that $\ell_I < \delta$ implies
that $O^{\sharp}$ exists, and, if $L^{\mu}$ exists, that also $O^{\dagger}$
exists. One aim of the paper is to give a modest improvement to this (see Theorem \ref{T3.2} and Cor. \ref{C3.4}).

\begin{definition}\label{D1.3} ({\em The $L^I$ validities}). We denote by $V_I$ the set of sentences of the logic $L^I$ true in all structures. 
\end{definition}

Notice that this notion of validity is over all such structures, not just those of a fixed signature. It is easy to see that the usual first order structure of arithmetic can be captured by a $L^I$ sentence, and from this it follows that the validities in $V_I$ are not arithmetically definable. In \cite{MR546441} it is shown that $V_I$ is neither $\Sigma^1_2$ nor $\Pi^1_2$ definable.

By an {\em Inner Model} we mean a transitive proper class model of the $\ZFC$ axioms, which thus contains $On$ - the class of all ordinals. Inner models which are of the form $L[E]$ ($E$ a coherent sequence of measures or extenders) are built under various assumtions concerning the size of inner models in the universe: as usual $O^\sharp$ denotes that there is no ``sharp for $L$" the constructible universe.
Further ``$\neg O^{\dagger}$'' abbreviates the assumption that there is no
sharp for any inner model with a single measurable cardinal; similarly let
``$\neg O^{\tmop{\pistol}}$'' (``not-Oh-pistol'') abbreviate the assumption that there is no sharp
for any inner model with a strong cardinal. 
We shall also use  ``$\neg O^{k}$'' (``not-Oh-kukri'', not ``not-Oh-kay'') for the assertion that there is no sharp for
 any inner model with a proper class of  measurable cardinals (see {\cite{We19}}).

Here $E$ is intended to be a coherent sequence of extenders; this can be in the sense of Jensen {\cf\!\!\cite{Je98}}, or as exposited in {\cite{zeman}}; this is constructed according
to such a sequence and is a very particular such model $K=L [ E^{K} ]$. 
The reader may also take this as the model built in \cite{St96}, although nothing in this article turns on the kind of indexing chosen for the extender $E$-sequence.
In fact here we shall
work with extenders suitable for building the {\tmem{core model}} $K$ below (or at) a
Woodin cardinal for which our official reference will be {\cite{JeSt13}}.

  The core models $K$ built under any of these anti-large cardinal assumptions,  may be very different.
  However for each of them (under the appropriate hypothesis) we have the Weak Covering Lemma (due in various models to Jensen and Mitchell, and below a Woodin cardinal, by Mitchell, Schimmerling, and Steel):\\

  { {\bf Weak Covering Lemma (\WCL(K))}} \cf\!\!\cite{zeman}, \cite{MiScSt97} {\em If $\gamma$ is a singular cardinal, then $\gamma$ is either singular or measurable in $K$. }\\

\setcounter{section}{3}

\setcounter{theorem}{0}

We show:
\begin{theorem}
  \label{test} \ \ Assume that there is no inner model of a Woodin cardinal, and that 
  there is a measurable cardinal.
Then $V_{I}$ is neither $\Sigma^{1}_{3}$ nor $\Pi^{1}_{3}$. \ 
\end{theorem}

Moreover: 
\setcounter{theorem}{11}
\begin{corollary}

Assume $V=L[E]$ but  there is no inner model of a Woodin cardinal. 
 Then $V_{I}$ is neither $\Sigma^{n}_{m}$ nor $\Pi^{n}_{m}$. 
\end{corollary}

\setcounter{theorem}{5}



In Section 3.2 we work under a further restricted notion of `smallness' of our models: we assume that no inner model (and so in particular no $L[E]$ model) has a measurable cardinal $\kappa$ which is strong up to some larger fixed point of the $\aleph$-function. We write this as  requiring that there is no inner model of ``$o(\kappa)=\aleph_{\mu}=\mu > \kappa$''. 
 We then can show that $K$ must have measurable cardinals:

\begin{corollary}  Assume $\ell_I<\delta$.
  Suppose  there is no inner model of ``$o(\kappa)=\aleph_{\mu}=\mu > \kappa$''. 
  Then the order type of the measurable
  cardinals of $K$
below the L\"owenheim number $\ell_I$ is  $\ell_I$.
\end{corollary}

We further have (Cor. \ref{C3.5}) a new example of Weak Covering.

\begin{corollary}  Assume $\ell_I<\delta$.
  Suppose  there is no inner model of ``$o(\kappa)=\aleph_{\mu}=\mu > \kappa$''.  Then the
  $\WCL(K)$ holds at the weakly inaccessible $\delta$ itself: $\delta^+=\delta^{+K}$.\\
\end{corollary}

\setcounter{section}{1}

 By $\ZFC^-$ we mean the Zermelo-Fraenkel axioms, taken with the axiom scheme of Collection, but with the power set axiom removed. By the {\em mouse ordering} we mean the relation $M \leq^{\ast} N$ which holds between set or proper classed-sized ``mice'' so that in the standard method of coiteration to comparable structures $(M_{\infty},N^{\infty})$ that $M_{\infty}$ is a (not necessarily proper) initial segment of the $N_{\infty}$ hierarchy. (See \cite{zeman}  again, Section 5.4.) We assume some slight familiarity with these notions and the construction there of $L[E]$ hierarchies with $E$ both a coherent sequence of measures, and also of extenders.

\section{Some results on core models and cardinals}


We first observe that any inner model
$M$ which has the same cardinals as $V$ will have the same core model as $V$ -
if those core models are thin: 


\begin{lemma}\label{Lem1}
  {\tmem{({\cite{We19}})}} Assume $\neg O^{k}$. Then $K (= K^{V})
  = K^{L [ \tmop{Card} ]}.$
\end{lemma}

  Consequently:
\begin{corollary} Assume $\neg O^{k}$ and that $M$ is an inner
  model with $\tmop{Card}  =  \tmop{Card}^{M}$, then\\ (i)  $K =K^{M}$;\\
  (ii) if additionally $V=L[E]$ then $M=V$.
\end{corollary}
\pf $K^{V}\, \sset \, L [ \tmop{Card} ]\, \sset \,M$.  Then the inductive construction of $K$ inside $M$ yielding $K^{M}$ in fact builds $K^{V}$. For (ii) $V=L[E]$ implies  
$K^{V}=V$; then $V\, \sset\, M$ by (i). \mb
\qed \\

Once $O^{k}$ exists then for such an $M$ as above, we may not have that $K^{M}=K$ 
but we shall nevertheless have that $K^{M}$ is {\tmem{universal}}. \
(Recall that a {\tmem{weasel}} $W$ is a class-sized mouse, and it is said to be {\em universal}
if in the comparison with any other mouse or weasel $P$ then $W$ absorbs $P$. 
That is if the comparison iteration of $( W,P )$ is to the models $( W_{\infty}
,P_{\infty} )$ then $P_{\infty}$ is an initial segment (not necessarily
proper) of $W_{\infty}$. \cf \cite{zeman} Sect. 6.3.)  A universal weasel is ``as good as'' being the full $K$ in many respects. Below a Woodin cardinal an inner model $M$ whose $K^{M}$ is universal will for example have the same reals as the true $K^{V}$.
If we strengthen the assumption to that of $\neg O^{\pistol}$ in Lemma \ref{Lem1.3}, then the universality of $K^{M}$ implies that it is a simple iterate of the true $K$. This is a result of R. Jensen and W. Mitchell, \cf\cite{St96}.

\begin{lemma}\label{Lem1.3}
  Suppose there is no inner model with a Woodin cardinal$\nosymbol$. Let $M$
  be any inner model with $\tmop{Card}^{M} = \tmop{Card}$, then $K^{M}$ is
  universal.
\end{lemma}

{\pf} Suppose this failed.  
We give here a standard  application of the Comparison Lemma ({\cite{MiSt91}} Theorem 7.1). Our supposition implies, with the proof of the Comparison Lemma, that
there is a cub class $C$ of points $i< \tmop{On}$ on the main branches  $b 
= [ 0, \infty ]_{\calu}$ and $  c =  [ 0, \infty ]_{\cal{T}}$
of the
iteration trees  $\mathcal{U}$, $\mathcal{T}$  resulting from the coiteration of
$N_{0} =K$ and $M_{0} =K^{M}$ \ to $( N_{\infty} ,M_{\infty} )$ with the following
properties, for $i<j \in C$, and $\pi^{N}_{i,j}$ ($\pi^{M}_{i,j} 
$ respectively) the iteration maps:\\

(i) $C \sset \, b,c $\,;\\

(ii) $\pi^{N}_{i,j} ( \kappa_{i} ) = \kappa_{j} =j$ ;\\

(iii) $\pi^{M}_{i,j} \m``\kappa_{j} \sset \kappa_{j}$.\\

Thus on the $K$ side some critical point is moved repeatedly by ultrapowers
out through the ordinals. However these ultrapower maps are continuous at the successor
of such critical points.  We thus also have that for some stage of the iteration $i_{0} \in C$ before which all truncations and all reductions in fine structural degree of the embeddings $\pi^{N}_{i,j}$ on the branch $c$ (if any) have occurred, that for later $j\in C\backslash i_{0}$:\\

(iv) $\pi^{N}_{i_{0},j} ( \kappa_{i_{0}}^{+N_{i}} ) = \kappa_{j}^{+N_{j}} \wedge  
\sup   \pi^{N}_{i_{0},j}   \m`` \kappa_{i_{0}}^{+N_{i_{0}}} =
\kappa_{j}^{+N_{j}}$.\\

(v) Thus there is some fixed $\tau \in
\tmop{Reg}^{V}$ so that \ $\tmop{cf}^{V} ( \kappa_{i}^{+N_{i}} ) = \tau$ for
$i_{0} \leq i \in C$.\\

\nod Let $\lambda >i_{0}$ be some larger limit  cardinal of $V$, with
$\kappa_{\lambda} = \lambda \in C$.\\ 

\nod{\tmem{Claim 1  $\l$ is inaccessible in $K^M$ (and all subsequent $M_{i}$ for $i\in C$).}\\

\pf Suppose $\l$ is singular in $K^{M}=M_{0}$. 
 If $f \in K^{M}$ is a
(1-1) increasing cofinal map with $f: \delta \imp \lambda$ for some $\delta <
\lambda$, then by our assumptions, by (ii) above as $\pi^{M}_{0, \lambda}
\m`` \kappa_{\lambda} \sset \kappa_{\lambda} =
\lambda$ we shall have $\tmop{ran} ( \pi^{M}_{0, \lambda} ( f ) ) \cap
\lambda$ \ a short increasing sequence witnessing the singularity of $\lambda$
in $M_{\lambda}$. But this is nonsense as in the coiteration $N_{\lambda}
\models$``$\lambda$ is measurable'' and $\power ( \lambda )^{N_{\lambda}} =
\power ( \lambda )^{M_{\lambda}}$. The argument for later $M_{i}$ is no different (or by normality of the iteration for $i\geq\l$). \qed Claim 1\\


\nod\tmem{Claim 2  $\pi^{M}_{0,\l}(\l) = \l$ and  $\sup\, \pi^{M}_{0,\l} \m``\l^{+M_{0}}= 
\pi^{M}_{0,\l}(\l^{+M_{0}}) = \l^{+ M_\l} $.}\\

\pf The first equality  holds by using (iii) that
 $\pi^{M}_{0,\l} \m``\l\sset \l$,
 after an argument by induction on $i<\l$  that $\pi^{M}_{0,i}(\l) = \l$; this uses the inaccessibility of $\l$ in the relevant models, and the usual argument as for measures, that if $i'$ is the $\calt$-successor of $i$, as $\pi^{M}_{i,i'}:M_{i}\imp Ult(M_{i},E_{\nu_{i}}) =
M_{i'}$ is a map from an extender of length $< \l$, then $\pi^{M}_{i,i'}$ fixes $\l$. This, using (iii), clearly holds into direct limits for $\calt$-limit $i'\leq_{\calt}
\l$.

For the second equality again by induction on $i<_{\calt}\l$  show that
$\pi^{M}_{0,i} \m``\l^{+M_{0}}$ is cofinal in $\l^{+M_{i}}$. Again this is the same argument as for measures and holds into direct limits for $\calt$-limit $i\leq_{\calt} \l$. The last equation is just by elementarity. \qed Claim 2\\

Now suppose $\l\in C$ had been chosen with $\tmop{cf}^{V} ( \lambda ) \neq \tau$ (where $\tau$ comes from (v)). By {\em Claim 2} $\tmop{cf}^{V}(\l^{+M_{0}}) =\tmop{cf}^{V}(\l^{+M_{\l}})$ and by the comparison process $\l^{+M_{\l}} < \lambda^{+}$
(because  $\l^{+M_{\l}}=  \l^{+N_{\l}}$ and the latter must be less than $ \lambda^{+}$).
Thus cofinality $\l^{+M_{\l}}$ has  cofinality $\tau$ in $V$. However by assumption on $\tmop{Card}^{M}$ we have $\l^{+}=\l^{+M}$ whilst $\l^{+K^{M}}<\l^{+M} $. As the \WCL\ holds in $M$ we cannot have that $\l$ is singular in $M$ and so, by \WCL\ again,
$(\tmop{cf}(\l^{+K^{M}}) = \l))^{M}$. Putting these facts together $\tmop{cf}^{V}(\l^{+K^{M}}) $ is now not equal to $\tau$. Contradiction! {\qed}






\begin{corollary}
  Suppose there is no inner model with a Woodin cardinal$\nosymbol$. Let $M$
  be any inner model so that:
  $$\{\l \in \tmop{Card} \mid \l^{+}=\l^{+M}\}\cap \tmop{Cof}_{\delta}
  $$
   is stationary for two different values of $\delta$,
    then $K^{M}$ is
  universal.
\end{corollary}
\pf If the conclusion failed then choose a value of $\delta$ which makes the given class stationary, but for a $\delta$ different from the $\tau$ of the last proof, which was the cofinality of the successor of the critical point used on the $N$ side, along the cub class $C$ contained in the main branches $b \cap c$.
\\ \mb \qed  \\

For our application later we remark that the comparison argument in the proof of Lemma \ref{Lem1.3} works for  sufficiently large $M$. 

\begin{corollary}\label{Cor1.5}
 Suppose there is no inner model with a Woodin cardinal$\nosymbol$. Let $M$
  be any transitive model of this statement together with a sufficiently large number of $\ZFC$ axioms, and with $\tmop{Card}^{M} = \tmop{Card}\cap M$, then $K^{M}$ is
{\em  weakly universal}: that is for any mouse $N$ with $On^{N}<On^{M}$ then $N<^{\ast}K^{M}$.
\end{corollary}
\pf We may assume $On^{M}$ is a strong limit cardinal and that sufficiently many axioms are true in $M$ to define $K^M$ and prove the $\tmop{\WCL}(K^M)$. Then taking $\theta$ as some singular cardinal below $On^{M}$, but with $On^{N}<\theta$ for a given $N$ as in the statement of the Corollary, we may assume the Weak Covering Lemma holds in $M$ of $K^{M}$, and in particular that $\theta^{+M} (= \theta^{+}) = \theta^{+K^{M}}$. Set $N_{0}=N$ and $P_{0}\dfs K^M\mid \theta^{+}$, and suppose
for a contradiction that $P_{0}\leq^{\ast} N_{0}$.  Then the comparison  to models $(N_{\theta^{+}},P_{\theta^{+}})$ requires $\theta^{+}$ steps, as indicated, with 
$P_{\theta^{+}}$ an initial segment of $N_{\theta^{+}}$, and there being no truncations of models or in the degree of ultrapowers taken on the main branch $[0,\theta^{+}]_{P_{0}}$ on the $P_{0}$-side. But now we obtain a contradiction, as $\theta^{+}$ must be a limit of critical points, and so inaccessibles in $N_{\theta^{+}}$, but cannot be so in $M_{\theta^{+}}$.
\qed\\

\section{Applications to strong logics}
\subsection{The complexity of the H\"artig logic validities}

The following is modelled on a proof that $V_{I}$ is neither $\Sigma^{1}_2$ nor $\Pi^{1}_{2}$ using Shoenfield Absoluteness of $\Sigma^{1}_{2}$ sentences of $L$ (\cf  \! \cite{MR546441}).

\begin{theorem}
  \label{test} \ \ Assume that there is no inner model of a Woodin cardinal, and that 
  there is a measurable cardinal.
Then $V_{I}$ is neither $\Sigma^{1}_{3}$ nor $\Pi^{1}_{3}$. \ 
\end{theorem}


{\pf}Let $\Upsilon$ be an $L^I$-sentence so that $\Upsilon^{M}$ holds iff
$M$ is (isomorphic to) a transitive model of $\tmop{\ZFC}_{N}$ (for some $N$,
$\tmop{\ZFC}_{N}$ an unspecified large number of axioms) with $\tmop{Card}^{M}
= \tmop{Card} \cap \tmop{On}^{M}$. Take such a transitive $M$ in which $\Upsilon^{M}$ holds. Let $\tau =
\tau ( M ) = \tmop{On}^{M}$. By the above Corollary \ref{Cor1.5}  $K^{M}$ is weakly universal, in particular, here for countable
mice (meaning it absorbs any mouse $N \in \tmop{HC}$). As any real of $K$ is
in some countable mouse, comparison of that mouse with $K^{M}$ shows that the
real is itself in $K^{M}$. Hence in particular $\re^{K} = \re^{K^{M}}$. Our
assumption on the existence of a measurable cardinal implies that $K$ is $\Sigma^{1}_{3}$-correct in $V$
({\cf}\cite{St96} Thm.7.9 - recall here that the second measurable cardinal $\Omega$ mentioned in this reference is only there to enable the construction of $K$; since the Jensen-Steel result of \cite{JeSt13} this upper measurable cardinal is redundant and a single measurable cardinal suffices). We thus have that
$\Sigma^{1}_{3}$-correctness in $V$ holds for $K^{M}$ and so $M$ too. Let
$\Phi ( n )$ define $S$, a complete $\Sigma^{1}_{3}$ set - which is perforce
not $\Pi^{1}_{3}$-definable. Then:
$$\Phi ( n )\Equi  \m``\Upsilon \imp  \Phi ( n ) ^{\tmop{HC}}\mq\in  V_{I}.$$
For, if $\Phi ( n )$ holds it will hold in $\tmop{HC}$ of any model $M$ with
$\Upsilon^{M}$ by the $\Sigma^{1}_{3}$-correctness we have just outlined. Hence
the quoted formula on the right hand side is $L^I$-valid, {\ie}it is in
$V_{I}$. Conversely, if $\neg \Phi ( n )$, then for a sufficiently large $\tau
\in \tmop{Card}$ we have: \ $(\Upsilon \wedge \neg \Phi ( n ))^{V_{\tau}}$, and thus the
right hand side fails. Hence $S$ is reducible to $V_{I}$ making the latter not
$\Pi^{1}_{3}$. By taking complements the same argument shows that $V_{I}$ is
not $\Sigma^{1}_{3}$.\qed\\

We shall improve this later in the case of $V$ being an $L[E]$ model at Corollary \ref{C3.9}.

\subsection{The L\"owenheim number for the H\"artig  logic, $\ell_{I}$}

In this subsection we show under a slightly more restrictive smallness assumption that the core model has measurable cardinals unbounded below  the L\"owenheim number $\ell_{I}$  for the H\"artig  logic.
\begin{theorem}\label{T3.2}
 (No inner model of ``$o(\kappa)=\aleph_{\mu}=\mu > \kappa$''.) Assume $\ell_I<\delta$. Then there are  measurable
  cardinals of $K$  unbounded
below the L\"owenheim number $\ell_I$. 
\end{theorem}

{\pf} Suppose not, and  the measurables of $K$ below $\ell_I$ are bounded by $\alpha_{0} <\ell_I$. We obtain a contradiction. Let $\Psi$ be a
sentence of $L^I$ that only has models of size at least
$\alpha_{0}$.  There is then a transitive model $M$ which is correct about cardinals, and in which there is an ordinal $\delta^{M}$ where $M\models \mbox{``}\delta^{M}\mbox{ is the least weakly inaccessible''}$. Further require that $M$ contains
a model $\mathfrak{A}$ of the sentence $\Psi$ which we may take as having cardinality in $M$ less than $\delta^{M}$. We may require $M$ to be a model of
$\tmop{\ZFC}_{n}$ a sufficiently large fragment of $\tmop{\ZFC}  $ that is
sufficient for the inductive construction of $K$ and to prove the Weak Covering Lemma
for it. Lastly require that in $K^{M}$ the measurables are bounded by the size of the model $\mathfrak{A}$, just as they are bounded by $\alpha_{0}$ in  $K$. Such a model can easily be found in $V$. 
By the definition of  $\ell_I$ we can then assume there is such an $M$ with  $\theta \dfs \tmop{On}^{M} \in
\tmop{Card}$, with $\alpha_{0} < \theta <\ell_I$.\\

We use the following lemma:

\begin{lemma}\label{L3.3}  (No inner model of ``$o(\kappa)=\aleph_{\mu}=\mu > \kappa$''.) 
Suppose (i) $N_{0},P_{0}$ are mice  with $P_{0}=L_{\gamma^{+}}[E]\models${\em ``There exists a largest cardinal''} for some $\gamma, \gamma^{+}\in Card$, and with  $N_{0}\in H_{\gamma^{+}}$. Suppose further that (ii) $P_{0}\models$``$\neg\exists \tau (o(\tau)\geq \gamma^{+})$''. Then $N_{0} <^{\ast}P_{0}$.  $P_{0}$ is thus universal for mice of smaller cardinality. 
\end{lemma}
\pf Suppose for a contradiction that $P_{0}\leq^{\ast} N_{0}$, and let $(N_{0},P_{0}) \longrightarrow (N_{
\theta},P_{\theta})$ be their coiteration, with $N_{0}$ iterating past $P_{0}$, and thus with $P_{\theta}$ an initial segment of $N_{\theta}$. We let the iteration maps be respectively $\pi^{N}_{i,j}$ and $\pi^{P}_{i,j}$ for $i\leq
j\leq\theta$. By this supposition there can be no truncation of any model $P_{i}$ on the $P$-side of the coiteration.  \\

{\em Claim}  $\pi^{P}_{{i,\gamma^{+}}}(\gamma^{+})=\gamma^{+}$.\\
Otherwise we should have for $i<\gamma^{+}$ that there is some $\lambda_{i}$ so that 
$\pi^{P}_{{i,\gamma^{+}}}(\lambda_{i})=\gamma^{+}$. By the usual arguments there is a c.u.b. $C\sset \gamma^{+}$ with $i<j\in C$ implying that $\pi^{P}_{{i,j}}(\lambda_{i})=\lambda_{j}$, and with $i=\lambda_{i}$ the critical point of the map  $\pi^{P}_{{i,j}}$. However $|N_{0}|<\gamma^{+}$ so likewise there is a c.u.b. $D\sset C$ unbounded in $\gamma^{+}$ with $\pi^{N}_{{i,j}}(\lambda_{i})=\lambda_{j}$ for  $i<j\in D$. By our  assumption (ii) on $P_{0}$, the comparison to $\gamma^{+}$ is the same as that for some $N_{i}$ and $P_{i}|\eta_{i}$ (where $\lambda_{i}$ is the critical point of $\pi^{P}_{i,i+1}$ and $\eta_{i}<\gamma^{+}$ is some ordinal with $\lambda_{i}^{+}\leq
\eta_{i}$). But these two mice are of cardinality $<\gamma^{+}$, and this coiteration is completed in less than $\gamma^{+}$ stages.  A contradiction and so the Claim holds.\\

However then $P_{\gamma^{+}}\models$``$\gamma^{+}$ is a successor cardinal'', whilst $N_{\gamma^{+}}$ believes that $\gamma^{+}$ is a limit of critical points $\lambda_{i}$, and so inaccessible. But this contradicts the comparison process.\\ \qed Lemma \ref{L3.3}


\begin{corollary}\label{C3.4}  (No inner model of ``$o(\kappa)=\aleph_{\mu}=\mu > \kappa$''.) $K_{\delta^{M}}$ $=^{\ast} K^{M}_{\delta^{M}}$. 

\end{corollary}

{\pf} Clearly $\mb^{\ast}\!\!\!\geq$ holds as $K_{\delta^{M}}$ is universal for all mice of cardinality $<\delta^{M}$. 
So suppose $K_{\delta^{M}}$ $\mb^{\ast}\!\!> K^{M}_{\delta^{M}}$. Again our smallness assumption  ensures there must be some $N_{0}$ a proper initial segment of $K_{\delta^{M}}$,  with both $|N_{0}|<\delta^{M}$ and $K^{M}_{\delta^{M}}\leq^{\ast}N_{0}$.
Let $\gamma > |N_{0}|$ be any singular cardinal with $\gamma<\delta^{M}$. The latter implies that $\gamma$ is singular in $M$, and thence, by $WCL$ in $M$, that $\gamma^{+}=\gamma^{+M}=\gamma^{+K^{M}}$. 
But now notice that it cannot be the case that for every such $\gamma$ satisfying the above that in $K^{M}$ there is some measurable $\tau_{\gamma}<\gamma$ with $K^{M}\models$``$o(\tau)\geq \gamma^{+}$'', since a regressive function argument would show then  that for some $\tau_{0}<\delta^{M}$ we'd have $K^{M}\models$``$o(\tau_{0})\geq \delta^{M}$''. Nor can it be the case that for all such $\gamma$ that $\gamma$ is measurable in $K^{M}$: for then in the comparison at stage $\gamma$ we should have that $\gamma^{+}=\gamma^{+P_{\gamma}}> \gamma^{+N_{\gamma}}$, and a truncation of $P_{\gamma}$ would be required on the $K^{M}$ side, which is impossible.

However then for some such $\gamma_{0}$ we have for $P_{0}=K^{M}\models$``$\gamma_{0} ^{+} = \gamma_{0}^{+K}$'',  to which we can apply the last lemma and the deduce that $N_{0}<^{\ast}K^{M}_{\gamma_{0}^{+}}$ - a contradiction.

\begin{corollary}  (No inner model of ``$o(\kappa)=\aleph_{\mu}=\mu > \kappa$''.) $K_\theta$ $=^{\ast} K^{M}$. 
\end{corollary}
\pf Similar. Note that we may assume that $M$ is a model of sufficiently many axioms so that $\theta = \aleph_{\theta}$. Hence  there is no $\tau<\theta$ with $(o(\tau)\geq \theta)^{K}$ to use as a smallness assumption once more.\\. \mb \qed \\



We  thus identify the coiteration of $( \nobracket K_{\delta^{M}}$, $K^{M}_{\delta^{M}} )   \nobracket$ to \ say $( \nobracket N_{\delta^{M}} ,
P_{\delta^{M}}$) as an initial part of the coiteration of $ \widetilde N_{0}\dfs K_{\theta} $ with $  \widetilde P_{0}\dfs K^{M}
)$ to some $( \widetilde N_{\infty} , \widetilde P_{\infty} )$, 
which lines up all the total measures to agreement and has the same sequence of indices. \\

\nod {\em Claim} At stage $\delta^{M}$ $\widetilde N_{\delta^{M}}$ must be truncated to some
$N^{\ast}_{\delta^{M}}$ to form an ultrapower to $\widetilde N_{\delta^{M} +1}$.\\

{\pf}Note at stage $\delta^{M}$ both models $\widetilde N_{\delta^{M}}$, $\widetilde P_{\delta^{M}}$ are of height the cardinal $\theta<\ell_{I}$.
 By the \WCL\ in
$V$ over $K$, $\delta^{M}$ being a singular cardinal, is either singular or
measurable in $K$. The latter fails by our assumption that such measurables below $\ell_{I}$ actually 
are all below $\alpha_{0}$, and our construction that enforces $\alpha_{0}\leq |\mathfrak{A}|<\delta^{M}$. Nevertheless $\delta^{M}$ is inaccessible in $M$
and thus so in $K^{M}$, and in the intermediate models $P_{\iota}$ for $\iota
\leq \delta^{M}$. However the iteration of $K_{\theta}$ to $\widetilde N_{\delta^{M}}$ preserves the singularity of $\delta^{M}$ in the models $\widetilde N_{\iota}$ (by induction on the stages $\iota \leq \delta^{M}$). Thus a truncation must be taken of $\widetilde N_{\delta^{M}}$ to
remove the subset of $\delta^{M}$ which is a witness to the singularity of $\delta^{M}$, if we are to have agreement between the final models. {\qed} Claim \\

However the result of this is that $K^{M}\leq^{\ast}N_{\delta^{M}+1}$, whilst the last corollary establishes that $K^{M}$, being $=^{\ast}$ equivalent to $K_{\theta}$, was universal for mice of cardinality $<\theta$. Contradiction! This  
establishes the theorem. \qed (Theorem)\\

\begin{corollary} \label{C3.6} Under the assumption of the theorem we in fact have: {\em  ``the order type of the $K$-measurables below $ \ell_{I}$ is $ \ell_{I}$''}.  
\end{corollary}
\pf
An assumption that the order type $\alpha_{0}<\ell_{I}$ can be used in the same way to obtain a contradiction. \qed \\

A further variant of the above yields:

\begin{corollary}\label{C3.5}  Assume $\ell_I<\delta$.
   Suppose  there is no inner model of ``$o(\kappa)=\aleph_{\mu}=\mu > \kappa$''. 
   Then:\\   (i)  the $\tmop{\WCL}(K)$ holds at $\delta$, {\ie} $\delta^{+} = \delta^{+K}$.\\
   (ii) there is $\delta\leq  \kappa < \delta^{++}$ which is measurable in $K$;\ \
\end{corollary}

{\pf}  Assume for a contradiction that one of these two conclusions fail. Run the above argument
but now using an $L^I$ sentence $\Psi_{1}$ stating:\\

 ``$\ex \dot{\delta} ( \dot{\delta} \nobracket$ {\em is the least weakly inaccessible
with either (a) $\dot{\delta}^{+} > \dot{\delta}^{+K}$  or  (b) no $\kappa$ with $\dot{\delta}\leq  \kappa < \dot{\delta}^{++}$
  measurable in $K)$}''.\\

{\nod}Then the model $M$ in the last argument, if it additionally is a model of $\Psi_{1}$, has cardinality some $\theta$ less than
$\ell_{I}$ with either, call it case (a), $\delta^{M+}= (\delta^M)^{+M}> (\delta^{M})^{+K^{M}}$, or otherwise, call it case   (b),  no $\kappa$ in $[\delta^M, (\delta^{M})^{++})$  measurable in $K^M$. Corollary \ref{C3.4} still holds as before.

Suppose case (a) could hold, then  by $\WCL(K)$ in $M$, we have: \\

\noindent {\em Either} $\delta^{M}$ is singular in $K=\widetilde N_{0}$ and so also in $\widetilde N_{\delta^{M}}$ and then we should have to do a truncation of $\widetilde N_{\delta^{M}}$  to remove the singularising sequence, as  $\delta^{M}$ inaccessible in $K^{M}$ implies that it is so in $\widetilde P_{\delta^{M}}$.

\noindent {\em Or } $\delta^{M}$ is measurable in $K$, with $ \delta^{M+} = ( \delta^{M} )^{+K}$. And it remains so in $\widetilde N_{\delta^{M}}$ with $( (\delta^{M}
)^{+})^{\widetilde N_{\delta^{M}}} = \delta^{M+}$, whilst $ (\delta^M)^{+\widetilde P_{\delta^{M}}}= (\delta^M)^{+\widetilde P_{0}}=(\delta^{M})^{+K^{M}} < (\delta^{M})^{+}$. So case (a) cannot hold and thus (i) is proven.\\

 Now
if case (b) were to hold then as $\delta^{M}$ is singular in $V$, by $\WCL(K)$ in $V$ then $\delta^{M}$ is singular or measurable in $K$.
If the former we have that $\widetilde N_{\delta^{M}}\models$``$\delta^{M}$ is singular'' whilst
 $\widetilde P_{\delta^{M}}\models$``$\delta^{M}$ is inaccessible''.
 Then we have to truncate  $N_{\delta^{M}}$ and iterate away a singularizing sequence which leads to a contradiction. If the latter, that is if $\delta^{M}$ is measurable in $K$, then  we should have that $N_{\delta^M}$ has $\delta^M$ as a measurable $K$-cardinal (as $\delta^M$ is not moved in the iteration of $N_0$ to $N_{\delta^M}$). However 
 then $N_{\delta^M}| {\delta^M}^{+}$ with its order zero measure, call it $E_{\delta^{M}}$, with critical point $\delta^{M}$, can neither be iterated out beyond $K^M_\theta$, (as above,  $K^M_\theta$ is universal for mice of cardinality less than $\theta$), nor, we shall argue, can it be absorbed as a measure on the $K^M$-side in some $P_{\ii}$ somewhere above $\delta^{M++}$, since the cofinality of the successor of its critical point $\delta^{M}$ is $\delta^{M+}$. The reason being that, if we set $\tau =\delta^{++}=\delta^{++M}$, then at the $\tau$'th stage of coiteration we should have that $\tau$ is measurable in $\widetilde N_{\tau}$. By agreement between the models we then have that $\tau^{+\widetilde P_{\tau}} < \tau^{+K^{M}}$. 
Moreover the map $\pi^{\widetilde P}_{0,\tau}\restriction \tau^{+K^{M}}$ is cofinal into $\tau^{+\widetilde P_{\tau}}$. But $M\models$``$cf( \tau^{+K^{M}})=\tau$''.
Hence $cf(\tau^{+\widetilde P_{\tau}})=\tau$
 whilst $cf( \tau^{+\widetilde N_{\tau}})= \delta^{+}<\delta^{++}=\tau$. Thus can again only be resolved by a truncation at stage $\tau$ on the $\widetilde N$ side. But this leads to a contradiction as before.
 
  Hence if this is to be  absorbed then it must happen at a stage, and so as a measure with critical point, below  $\delta^{M++}$. 
But then
 a truncation of $\widetilde N_{\delta^{M}}$ is again required and this is a contradiction as before. 
 {\qed}\\

From these conclusions one might expect that it is indeed $\delta$ itself which is measurable in $K$, but the proof falls just short of that.
However it seems we can push this argument further:

\begin{theorem}
  Assume that \ $\ell_I < \delta$.   Suppose  there is no inner model of ``$o(\kappa)=\aleph_{\mu}=\mu > \kappa$''
 Then for every weakly inaccessible
  cardinal $\lambda$ there is $\lambda\leq
  \kappa <\lambda^{++}$ which is measurable in $K$ and moreover $\lambda^{+} = \lambda^{+K}$.
\end{theorem}

{\pf} We assume the hypothesis, and the smallness supposition. The argument of the last Corollary used no special properties of
$\delta$, for example that it was the least weakly inaccessible. So we may
repeat this to obtain a contradiction for any weakly inaccessible cardinal
$\lambda$ for which the conclusions supposedly fail. {\qed}\\

Thus a proper class of weakly inaccessibles under this hypothesis yields a proper class of measurables in $K$. (We should add here that it is unknown if this hypothesis is consistent: the argument of \cite{MaVa11} only yields a $ZFC$ model with a single weakly inaccessible $\delta> \ell_I$.)

\subsection{When the power set operation is {\boldmath$\Sigma_{1}(\tmop{Card})$} definable}

\begin{lemma}\label{L2.3}
  Suppose $V=L [ E ]$ is a model of {\em ``There is no inner model of a Woodin cardinal.''}
  Then the relation ``$x=\mathcal{P} ( y )$'' \ is $\Sigma_{1} ( \tmop{Card} )$ where
  $\tmop{Card}$ is the predicate true of just the infinite cardinals. 
\end{lemma}
 \pf It suffices to prove the Lemma for $y\in Card$.
 So, let $\alpha \geq \omega_{2}$, $\alpha \in \tmop{Card}$. 
 We use the following folk-lore style lemma, which is taken from \cite{ClSch12}, Lemma 2.1. This ensures that a premouse is a mouse, \ie is $\omega_1 +1$-iterable, if it is seen so by a sufficiently closed transitive model.  The clause below that requires there be no definably Woodin cardinals, means that any Woodin cardinal of $M$ is either definably collapsed over $M$ or there is a counterexample to its Woodiness, again definable over $M$.

 \begin{lemma}\label{CS}
 Suppose there is no inner model of a Woodin cardinal. Let $U$ be an uncountable transitive model of this statement, together with $\ZFC^-$. Let $M \in
 U$ be a premouse, which has no definably Woodin cardinals. Then $M$ is a mouse if and only if $(M \mbox{ is a mouse })^U$.
 \end{lemma}
We shall also use version of the Lemma 3.5 of \cite{GiShSch06}:

\begin{lemma}Suppose there is no inner model of a Woodin cardinal. Let $K$ denote the core model. Let $\kappa \geq \aleph_{2}$ be a cardinal of K, and let $N \triangleright K| | \kappa$ be an iterable premouse such that $\rho^{\omega}_{N} \leq \kappa$, and $N$ is sound above $\kappa$. Then $N\triangleleft K$, \ie $N$ is an initial segment of $K$.
\end{lemma}
(We have suppressed mention of the larger measurable cardinal $\Omega$ in the Lemma as this is now redundant.
Also the lemma there is stated in a much wider form, and we have taken the base case (``n=0") only. As written in \cite{GiShSch06}, this would require  that $V$ be closed under $\sharp$'s, but in fact this is not needed for the base case.)\\

(1) Let $N= \langle J^{E^{N}}_{\delta} ,E^{N} \rangle$ be a 
premouse with $\delta \in \tmop{LimCard}$ and $\tmop{Card}^{N} = \tmop{Card} \cap
\tmop{On}^{N}$. Suppose $N\models$``{\em Every initial segment $\pa{J^{E^{N}}_{\alpha} ,E^{N}, E^N_\alpha}$ is iterable''. }
Then $N = K_\delta = L_{\delta}[E]$.\\

\pf The point is that iterability for such  
$N$ has been turned into a first order statement. First note that  every set $x\in N$ is an element of a transitive $\ZFC^-$ model; secondly that $N$ is a model of {\em ``there is no inner model of a Woodin cardinal''}  (suppose $J_\delta^{\bar E}$ were an inner model of $N$ with a Woodin cardinal. As $\delta \in \tmop{Card}$ we should have that $L[\bar E]$ was an inner. model of a Woodin).  Thirdly, $N$ itself can have no Woodin cardinals, for the same reason. Hence $N$ has no definably Woodin cardinals. Thus we can just apply  Lemma \ref{CS} to see that all initial segments  $\pa{J^{E^{N}}_{\alpha} ,E^{N}, E^N_\alpha}$ are mice.

Now consider the putative comparison of $N= N_0$ with $M_0 =J_\delta^E$. We see immediately that $N|\omega_2 = J^E_{\omega_2}$ since there are no measurable cardinals in these structures to coiterate, and any comparison could not truncate on both sides of the coiteration to make them, and a truncation on one side only results in a host of inaccessible cardinals which are not there in the other side. Hence $L_{\omega_2}[E]= K_{\omega_2} = N|\omega_2$. But now by induction on $K$-cardinals $\omega_2<\kappa \leq \delta$ (which are just the $V=L[E]$ cardinals) we show that $K|\kappa$ is an initial segment of $N$. (This is because the premouse $N$ has the same cardinals as $L_\delta[E]$ and thus any $\kappa^{+N}$ is the supremum
 of ordinals $\alpha$ of initial segments $J_\alpha^{E^N}$ with projectum dropping to $\kappa$. Thus such initial segments are trivially in $K$ and so too on the $K$-sequence $E$. And again $\kappa^{+N}= \kappa^+$.)  Thus $N=K_\delta$. \qed (1)\\

\nod We just note that this will finish the Lemma.  Setting $\Psi$ to be ``Every initial segment is iterable'':

$$\begin{array}{lcl}
 y=\mathcal{P} ( \alpha ) &\Leftrightarrow&  
 \exists N ( N  \mbox{  is a premouse, }
\Psi^N \wedge
 \tmop{On}^{N}
\in \tmop{LimCard} \wedge \\
&&\tmop{Card}^{N} = \tmop{Card} \cap
\tmop{On}^{N} \wedge ( y=\mathcal{P} ( \alpha ) )^{N} ).
\end{array}$$
{\noindent}as being a premouse is a $\Delta_{1}$ notion, being defined by just
first order properties over the structure, and the rest considered as $\Sigma_1(\tmop{Card})$. \qed (Lemma \ref{L2.3})\\

\begin{corollary}\label{C3.9}

Assume $V=L[E]$ but  there is no inner model of a Woodin cardinal. Then $V_{I}$ is neither $\Sigma^{n}_{m}$ nor $\Pi^{n}_{m}$. 
\end{corollary}
\pf  Now let $\Upsilon_1$ be an $ L^I$-sentence so that $\Upsilon_1^{M}$ holds iff
$M$ is (isomorphic to) a transitive model of $\tmop{\ZFC}_{N}$ (for some $N$ large),
with $\tmop{Card}^{M}
= \tmop{Card} \cap \tmop{On}^{M}$,
and is also a sound premouse which thinks all its initial segments are iterable mice.
Take such a transitive $M$ in which $\Upsilon^{M}$ holds.


 Then such an $M$ is correct about power sets, indeed has domain $V_\delta = L_\delta[E]$, in particular is obviously correct about finitely iterated power sets of $\nat$ by the last Lemma. Then use the template of Theorem \ref{test}. \qed \\

Obviously this can be extended to show the undefinability of $V_{I}$ over much higher types.
\section{Open Questions}

\nod {\tmem{Question 1}}: Can we improve the lower bound in Theorem \ref{T3.2}? The consistency of $\ell_I$ less than the first weakly inaccessible is obtained from a supercompact cardinal. There is thus a wide gap here.
\\

\nod {\tmem{Question 2}}: How much of the above works for the logic $L(I,Q^{ec})$ \cf \cite{MaVa11}. We may be able to reflect down below the first weakly Mahlo, but can we get further measurables in $K$ as a result?
\\

\nod {\tmem{Question 3}}: How large can a cardinal be and still be consistent with the statement ``$x=\mathcal{P} ( y )$ is  $\Sigma_1(Card)$"?
\\

\nod{\tmem{Question 4}}: Is it consistent relative to the existence of a supercompact cardinal,  that there is no proper inner model $M$ with $Card^M=Card^V$?

\small
\bibliographystyle{plain}
\bibliography{settheory10t,patras}

\begin{thebibliography}{10}

\bibitem{ClSch12}
B.~Claverie and R-D. Schindler.
\newblock {Woodin's} axiom $(\ast)$, bounded forcing axioms, and precipitous
  ideals on $\omega_1$.
\newblock {\em J. Symbolic Logic}, 77(2):475--498, 2012.

\bibitem{MR0406772}
Solomon Feferman.
\newblock Applications of many-sorted interpolation theorems.
\newblock In {\em Proceedings of the {T}arski {S}ymposium ({P}roc. {S}ympos.
  {P}ure {M}ath., {V}ol. {XXV}, {U}niv. of {C}alifornia, {B}erkeley, {C}alif.,
  1971)}, pages 205--223, 1974.

\bibitem{GiShSch06}
M.~Gitik, S.~Shelah, and R-D Schindler.
\newblock Pcf theory and {Woodin} cardinals.
\newblock In {\em Logic Ciolloquium '02}, volume~27 of {\em Lecture Notes in
  Logic}, pages 172--205. Association for Symbolic Logic, 2006.

\bibitem{Je98}
R.~B. Jensen.
\newblock A new fine structure for higher core models.
\newblock {\em Circulated manuscript}, Berlin, 1997.

\bibitem{JeSt13}
R.B. Jensen and J.~R. Steel.
\newblock {$K$} without the measurable.
\newblock {\em Journal of Symbolic Logic}, 78(3):708--734, Sep. 2013.

\bibitem{MR244012}
Per Lindstr\"{o}m.
\newblock First order predicate logic with generalized quantifiers.
\newblock {\em Theoria}, 32:186--195, 1966.

\bibitem{MaVa11}
M.~Magidor and J.~V{\"{a}}{\"{a}}n{\"{a}}nen.
\newblock On {L\"owenheim}-{Skolem}-{Tarski} numbers for extensions of first
  order logic.
\newblock {\em JML}, 11(1):87--113, 2011.

\bibitem{MR457146}
J.~A. Makowsky, Saharon Shelah, and Jonathan Stavi.
\newblock {$\Delta$}-logics and generalized quantifiers.
\newblock {\em Ann. Math. Logic}, 10(2):155--192, 1976.

\bibitem{MiScSt97}
W.~J. Mitchell, E.~Schimmerling, and J.R. Steel.
\newblock The covering lemma up to a {Woodin} cardinal.
\newblock {\em Annals of Pure and Applied Logic}, 84:219--255, 1997.

\bibitem{MiSt91}
W.~J. Mitchell and J.~R. Steel.
\newblock {\em Fine Structure for Iteration Trees}, volume~3 of {\em Lecture
  Notes in Logic}.
\newblock Springer-Verlag, Berlin, 1991.

\bibitem{MR1928384}
J.~Stavi and J.~V\"a\"an\"anen.
\newblock Reflection principles for the continuum.
\newblock In {\em Logic and algebra}, volume 302 of {\em Contemp. Math.}, pages
  59--84. Amer. Math. Soc., Providence, RI, 2002.

\bibitem{St96}
J.~R. Steel.
\newblock {\em The Core Model iterability problem}, volume~8 of {\em Lecture
  Notes in Mathematical Logic}.
\newblock Springer, 1996.

\bibitem{MR510714}
Jouko V\"a\"an\"anen.
\newblock Two axioms of set theory with applications to logic.
\newblock {\em Ann. Acad. Sci. Fenn. Ser. A I Math. Dissertationes}, (20):19,
  1978.

\bibitem{MR546441}
Jouko V\"{a}\"{a}n\"{a}nen.
\newblock Remarks on free quantifier variables.
\newblock In {\em Essays on mathematical and philosophical logic ({P}roc.
  {F}ourth {S}candinavian {L}ogic {S}ympos. and {F}irst {S}oviet-{F}innish
  {L}ogic {C}onf., {J}yv\"{a}skyl\"{a}, 1976)}, volume 122 of {\em Synthese
  Library}, pages 267--272. Reidel, Dordrecht-Boston, Mass., 1979.

\bibitem{MR654791}
Jouko V\"{a}\"{a}n\"{a}nen.
\newblock Abstract logic and set theory. {II}. {L}arge cardinals.
\newblock {\em J. Symbolic Logic}, 47(2):335--346, 1982.

\bibitem{MR694269}
Jouko V\"a\"an\"anen.
\newblock Generalized quantifiers in models of set theory.
\newblock In {\em Patras {L}ogic {S}ymposion ({P}atras, 1980)}, volume 109 of
  {\em Stud. Logic Foundations Math.}, pages 359--371. North-Holland,
  Amsterdam-New York, 1982.

\bibitem{We19}
P.D. Welch.
\newblock Closed and unbounded classes and the {H\"artig} quantifier model.
\newblock {\em J. Symbolic Logic (to appear)}, page~18.

\bibitem{zeman}
M.~Zeman.
\newblock {\em Inner models and large cardinals}, volume~5 of {\em de Gruyter
  Series in Logic and its Applications}.
\newblock Walter de Gruyter \& Co., Berlin, 2002.

\end{thebibliography}

\end{document}